\documentclass{amsart}
\usepackage{amssymb, latexsym}

\newcommand{\REFER}{\ref}

\newcommand{\R}{{\mathbb R}}
\newcommand{\C}{{\mathbb C}}

\newcommand{\CC}{{\overline \C}}
\newcommand{\CR}{{\overline \R}}

\newcommand{\nbhd}{{neighborhood }}

\newtheorem{conj}{Conjecture}  
\newtheorem{lemma}{Lemma}
\newtheorem{property}{Property}
\newtheorem{theorem}{Theorem}

\theoremstyle{definition}
\newtheorem{definition}{Definition}
\theoremstyle{remark}
\newtheorem{remark}{Remark}

\newtheorem{example}{Example}

\numberwithin{equation}{section} \numberwithin{theorem}{section}
\numberwithin{definition}{section} \numberwithin{remark}{section}
\numberwithin{example}{section} \numberwithin{lemma}{section}
\numberwithin{property}{section}
\numberwithin{proposition}{section} \numberwithin{claim}{section}
\numberwithin{othertheorem}{section} \numberwithin{conj}{section}
\numberwithin{corollary}{section}

\begin{document}

\title{Some counterexamples in dynamics of rational semigroups}

\author{Rich Stankewitz}
\address[Rich Stankewitz]
    {Department of Mathematical Sciences\\
    Ball State University\\
    Muncie, IN 47306\\
    USA}
    \email{rstankewitz@bsu.edu}

\author{Toshiyuki Sugawa}
\address[Toshiyuki Sugawa]
{Department of Mathematics\\
Graduate School of Science\\
Hiroshima University\\
Higashi-Hiroshima, 739-8526 Japan}
\email{sugawa@math.sci.hiroshima-u.ac.jp}

\author{Hiroki Sumi}
\address[Hiroki Sumi]{Department of Mathematics\\
    Tokyo Institute of Technology 2-12-1\\
    O-Okayama, Meguro-ku\\
    Tokyo, 152-8551, Japan}
    \email{sumi@math.titech.ac.jp}

%

\thanks{
The second author was partially supported by the Ministry of
Education, Grant-in-Aid for Encouragement of Young Scientists,
11740088. }


\begin{abstract}
We give an example of two rational functions with non-equal Julia
sets that generate a rational semigroup whose completely invariant
Julia set is a closed line segment.  We also give an example of
polynomials with unequal Julia sets that generate a non nearly
Abelian polynomial semigroup with the property that the Julia set
of one generator is equal to the Julia set of the semigroup. These
examples show that certain conjectures in the field of dynamics of
rational semigroups do not hold as stated and therefore require
the allowance of certain exceptional cases.
\end{abstract}

\maketitle

\section{Introduction}

In~\cite{HM1}, Hinkkanen and Martin develop a theory of dynamics
of rational semigroups as a generalization of the classical theory
of the dynamics of the iteration of a rational function defined on
the Riemann sphere $\CC$.  In that paper and in subsequent
communications, they put forth several conjectures, some of which
will be addressed here. In particular, we provide counterexamples
to Conjectures \ref{conj1}, \ref{conj2} and \ref{Conj:Jf=Jg}.  In
light of these examples the conjectures are then suitably modified
and as such remain open questions.  We begin by developing the
necessary background to state these questions.

In what follows all notions of convergence will be with respect to
the spherical metric on $\CC.$ A rational
semigroup $G$ is a semigroup of rational functions of degree at
least two defined on $\CC$ with the semigroup
operation being functional composition. (One may wish to allow
some or all of the maps in $G$ to be M\"obius, for example, when
one is considering Kleinian groups as in \cite{RS3}, but since the
examples constructed here all contain maps of degree two or more,
we will use our simplified definition to avoid any technical
complications which are not pertinent to this paper.) When a
semigroup $G$ is generated by the functions $\{f_1, f_2, \dots,
f_n, \dots\}$, we write this as
\begin{equation}
G=\langle f_1, f_2,\dots, f_n, \ldots \rangle.\notag
\end{equation}

On p. 360 of~\cite{HM1}, the definitions of the set of normality,
often called the Fatou set, and the Julia set of a rational
semigroup are as follows:

\begin{definition} \label{N(G),J(G)}
For a rational semigroup $G$ we define the set of normality of
$G$, $N(G)$, by
$$N(G)=\{z \in \CC:\exists \text{ a \nbhd of } z \text{ on
which } G \text{ is a normal family} \}$$ and define the Julia set
of $G$, $J(G)$, by
$$J(G)=\CC \setminus N(G).$$
\end{definition}

Clearly from these definitions we see that $N(G)$ is an open set
and therefore its complement $J(G)$ is a compact set. These
definitions generalize the case of iteration of a single rational
function and we write $N(\langle h \rangle )=N_h$ and $J(\langle h
\rangle )=J_h$. Note that $J(G)$ contains the Julia set of each
element of $G$.  For research on (semi-)hyperbolicity and
Hausdorff dimension of Julia sets of rational semigroups,
see~\cite{Su1,Su3,Su2,Su4}.

\begin{definition} If $h$ is a map of a set $Y$ into
itself, a subset $X$ of $Y$ is:
\begin{align}
  i)& \,forward\,\,invariant \text{ under } h \text{ if } h(X) \subset
  X;\notag\\ ii)& \,backward\,\,invariant \text{ under } h \text{ if }
  h^{-1}(X) \subset X;\notag\\ iii)& \,completely\,\,invariant \text{
  under } h \text{ if } h(X) \subset X \text{ and } h^{-1}(X) \subset
  X.\notag
\end{align}

\end{definition}

It is well known that for a rational function $h$, the set of
normality of $h$ and the Julia set of $h$ are completely invariant
under $h$ (see ~\cite{Be}, p.~54), i.e.,
\begin{equation}\label{cinv}
h(N_h)=N_h=h^{-1}(N_h) \text{ and } h(J_h)=J_h=h^{-1}(J_h).
\end{equation}

In fact, the following property holds.

\begin{property}\label{Jh2}
For a rational map $h$ of degree at least two the set $J_h$ is the
smallest closed completely invariant (under $h$) set which
contains three or more points (see ~\cite{Be}, p.~67).
\end{property}

From Definition~\ref{N(G),J(G)}, it follows that $N(G)$ is forward
invariant under each element of $G$ and, thus, $J(G)$ is backward
invariant under each element of $G$ (see ~\cite{HM1}, p.~360). The
sets $N(G)$ and $J(G)$ are, however, not necessarily completely
invariant under the elements of $G$.  This is in contrast to the
case of single function dynamics as noted in (\ref{cinv}).
However, one could generalize the classical notion of the
Julia set of a single function in such a way as to force the Julia
set of a rational semigroup to be completely invariant under each
element of the semigroup. Thus, we give the following definition.

\begin{definition}\label{Edef}
 For a rational semigroup $G$ we define the \textbf{completely invariant}
 Julia set of $G$
$$E(G)=\bigcap\{S:S \text{ is closed, completely invariant under
each }g \in G, \#(S)\geq3 \}$$ where $\#(S)$ denotes the
cardinality of $S$.
\end{definition}

We note that $E(G)$ exists, is closed, is completely invariant
under each element of $G$ and contains the Julia set of each
element of $G$ by Property~\REFER{Jh2}.

\begin{definition}\label{Wdef}
For a rational semigroup $G$ we define the \textbf{completely
invariant} set of normality of $G$, $W(G)$, to be the complement
of $E(G)$, i.e.,
$$W(G)=\CC \setminus E(G).$$
\end{definition}

Note that $W(G)$ is open and it is also completely invariant under
each element of $G$.


We state the following conjectures which are due to A. Hinkkanen
and G. Martin (see~\cite{RS2}).

\begin{conj}\label{conj1}
If $G$ is a rational semigroup which contains two maps $f$ and $g$
such that $J_f \neq J_g$ and $E(G) \neq \CC$, then $W(G)$ has
exactly two components, each of which is simply connected, and
$E(G)$ is equal to the boundary of each of these components.
\end{conj}

\begin{conj}\label{conj2}
If $G$ is a rational semigroup which contains two maps $f$ and $g$
such that $J_f \neq J_g$ and $E(G) \neq \CC$, then $E(G)$ is a
simple closed curve in $\CC$.
\end{conj}

In section~\ref{counterex1} we give a method for constructing
functions (as well as providing concrete functions) whose Julia
sets are unequal, but which generate a semigroup whose completely
invariant Julia set is a line segment.  Hence the above
conjectures do not hold.  But since the only completely invariant
Julia sets of rational semigroups which are known at this time
(when the semigroup contains two maps with unequal Julia sets) are
$\CC$ (see~\cite{RS1} and \cite{RS2}) or sets which are M\"obius
equivalent to a line segment or circle, the authors put forth the
following conjecture, which is currently unresolved.

\begin{conj}\label{conj1'}
If $G$ is a rational semigroup which contains two maps $f$ and $g$
such that $J_f \neq J_g$ and $E(G)$ is not the whole Riemann
sphere, then $E(G)$ is M\"obius equivalent to a line segment or a
circle.
\end{conj}

\begin{remark}
We briefly explain some evidence that compels us to pose
Conjecture~\ref{conj1'} in this way.  Our example of a rational
semigroup $G$ with $E(G)$ being a line segment is rigid since $G$
contains a Tchebycheff polynomial, which is known to be
postcritically finite (and hence, rigid). On the other hand, an
example of a rational semigroup $G$ with $E(G)$ being a (unit)
circle generated by rational functions $f_1,\dots,f_n$ with
non-equal Julia sets is easily constructed by choosing finite
Blaschke products as the $f_j$'s. However, it seems difficult to
quasiconformally deform $f_1,\dots,f_n$ simultaneously so that the
completely invariant Julia set of the resulting rational semigroup
is not a circle.
\end{remark}

In section~\ref{counterNA} we provide a counterexample to the
following conjecture also due to Hinkkanen and Martin~\cite{HMpc}.
\begin{conj} \label{Conj:Jf=Jg}
         Let $G$ be a polynomial semigroup such that $J_h = J(G)$ for some
         $h \in G$.  Then $J_f = J_g$ for all $f,g \in G$ (and hence $G$ is
         nearly abelian by Theorem~\ref{C:polysameJ}).
\end{conj}

In our counterexample $J(G)$ is a closed line segment.  Since no
other types of counterexamples are known, we modify this
conjecture as follows and note that it remains unresolved.

\begin{conj} \label{Conj:Jf=Jg'}
         Let $G$ be a polynomial semigroup such that $J_h = J(G)$ for some
         $h \in G$ where $J(G)$ is not a line segment.
         Then $J_f = J_g$ for all $f,g \in G$ (and hence $G$ is
         nearly abelian by Corollary~\ref{C:polysameJ}).
\end{conj}

\section{Counterexamples to Conjectures~\ref{conj1} and~\ref{conj2}}
\label{counterex1}


We begin this section with some notation and lemmas.  Let
$\phi(z)=\frac{z^2-1}{z^2+1}$ and denote the upper half plane as
$U=\{z:\Im z>0\}$. Then $\phi$ maps $U$ one-to-one onto
$\Omega=\CC \setminus [-1,1]$ and $\phi$ maps $\CR$ two-to-one
onto $I=[-1,1]$. We call a map $f$ \textit{odd} if $f(-z)=-f(z)$
and we call a map $f$ \textit{even} if $f(-z)=f(z)$.

\begin{lemma}\label{form}
A function $f$ is an odd rational map such that $f(U)=U$ if and
only if it has the form

\begin{equation}\label{eqnform}
f(z)=az - \frac{b}{z} - \sum_{j=1}^N \frac{B_j z}{z^2-A_j}
\end{equation}

 where
$a, b, A_j, B_j \geq 0.$
\end{lemma}

\begin{proof}
Let $f$ be an odd rational map such that $f(U)=U$. Then any
preimage of infinity must be real (else there would exist a
preimage of infinity in $U$) and simple (else there would be
points in $U$ that map outside of $U$). Again, since $f(U)=U$, it
follows that $f$ must be of the form $az-\frac{b}{z} -
\sum_{j=1}^k \frac{c_j}{z-a_j}$ where $a, b, c_j \geq 0$ and $a_j
\in \R$. Since $f(\CR)=\CR$ and $f(-z)=-f(z)$ we conclude that the
poles other than the one which might possibly exist at the origin
must come in pairs of real numbers symmetric about the origin.
Hence $f(z)=az-\frac{b}{z} - \sum_{j=1}^N \frac{b_j}{z-a_j} -
\sum_{j=1}^N \frac{b_j}{z+a_j}$ where $b_j >0$, which can be
algebraically reduced to~\eqref{eqnform}.

Let $f$ be a map of the form~\eqref{eqnform}.  Hence $f$ is odd,
rational, maps $U$ into $U$ (since each term in the sum does), and
maps $\CR$ into $\CR$ (since the coefficients are all real). From
this it easily follows that $f(U)=U$.
\end{proof}

\begin{lemma}\label{fevenorodd}
Let $f$ be a rational map.  Then $[f(z)]^2$ is even if and only if
$f$ is even or odd.
\end{lemma}

\begin{proof}
Suppose $[f(z)]^2=[f(-z)]^2$.  Then an analytic square root
 of $[f(-z)]^2$ (defined locally away from the zeroes and poles of $f$)
is either $f(z)$ or $-f(z)$.  The identity theorem can then be
used to show that $f(-z)$ is either $f(z)$ or $-f(z)$ globally,
i.e., $f$ is either even or odd.

The reverse implication is immediate.
\end{proof}

\begin{lemma}\label{f2even}
Let $f$ be a rational map.  Then $\phi \circ f$ is even if and
only if $[f(z)]^2$ is even.
\end{lemma}

\begin{proof}
Since $\phi(z)=\psi(z^2)$ for $\psi(z)=\frac{z-1}{z+1}$ we see
that $[f(z)]^2=\psi^{-1}\circ \phi(f(z))$ and $\psi^{-1}\circ
\phi(f(-z))=[f(-z)]^2.$  The lemma easily follows.
\end{proof}

\begin{lemma}\label{geven}
If $g$ is an even rational function, then $g(z)=h(z^2)$ for some
rational map $h$.
\end{lemma}

\begin{proof}
For $z\neq 0$ or $\infty$ we define $h(z)=g(\pm \sqrt{z})$ and
note that $h$ is well defined (regardless of the branch of the square root taken)
since $g$ is even.  Since $h$ is
analytic on $\C \setminus \{0\}$ and can be extended in the
obvious way to be continuous on $\CC$, $h$ is rational and
satisfies $h(z^2)=g(z)$.
\end{proof}

\begin{lemma}\label{tildef}
Let $f$ be a rational map such that $f(U)=U$.  Then there exists a
rational map $\tilde{f}$ such that $\phi \circ f = \tilde{f} \circ
\phi$ if and only if $f$ is odd (and therefore of the form in
Lemma~\ref{form}).
\end{lemma}

\begin{proof}
Let $f$ be odd.  Since $f(z)=-f(-z)$ we see that $[f(z)]^2$ is an
even rational function and therefore by Lemma~\ref{geven}
$[f(z)]^2=h(z^2)$ for some rational map $h$.  Define
$\tilde{f}(z)=\frac{h(\frac{1+z}{1-z})-1}{h(\frac{1+z}{1-z})+1}$
(hence $\tilde{f}$ is a rational map as it is a composition of
rational maps).  Let $w=\phi(z)=\frac{z^2-1}{z^2+1}$ and note that
$z^2=\frac{1+w}{1-w}$.   Hence $(\tilde{f} \circ
\phi)(z)=\tilde{f}(w)=\frac{h(\frac{1+w}{1-w})-1}{h(\frac{1+w}{1-w})+1}
=\frac{[f(z)]^2 -1}{[f(z)]^2+1}= (\phi \circ f) (z)$.

Suppose there exists a rational map $\tilde{f}$ such that $\phi
\circ f = \tilde{f} \circ \phi$.  Then $\tilde{f}\circ \phi$ is
even since $\phi$ is even.  The semi-conjugacy implies $\phi \circ
f$ is also even, which by Lemmas~\ref{f2even} and \ref{fevenorodd}
gives that $f$ is either even or odd.  If $f$ were even, then
$f(\CC \setminus U)=f(\overline{U})=\overline{U}$ and the preimage
of the lower half plane would be empty.  This contradicts the fact
that the image of $\CC$ under a rational map is always $\CC$.
Hence we conclude that $f$ must be odd.
\end{proof}
\begin{lemma}\label{lift}
If $\tilde{f}$ is a rational map such that
$\Omega=\tilde{f}^{-1}(\Omega)$, then there exists an odd rational
map $f$ such that $U=f^{-1}(U)$ and $\phi \circ f = \tilde{f}
\circ \phi$.
\end{lemma}

\begin{proof}
Let $h$ denote the branch of the inverse of $\phi$ which maps
$\Omega$ onto $U$.  Then $f= h \circ \tilde{f} \circ \phi$ maps
$U$ onto $U$ properly and is therefore a rational map (Blaschke
product of the upper half plane).  Clearly, $\phi \circ f =
\tilde{f} \circ \phi$ on $U$ and so by the identity Theorem this
semi-conjugacy holds on all of $\CC$.  By Lemma~\ref{tildef} $f$
is odd.
\end{proof}

\begin{remark}
Lemmas~\ref{tildef} and \ref{lift} classify those rational
functions that can be semi-conjugated by $\phi$.
\end{remark}

\begin{lemma}\label{semiconj}
For rational semigroups $G=\langle g_j:  j \in \mathcal{I}
\rangle$ and $H=\langle h_j: j \in \mathcal{I} \rangle$ where
there exists a rational function $k$ satisfying the semi-conjugacy
relation $k \circ h_j=g_j \circ k$ for each $j \in \mathcal{I}$,
we have $J(G)=k(J(H))$ and $N(G)=k(N(H))$.
\end{lemma}

\begin{proof}
We first note that the semi-conjugacy relation on the generators
translates to a semi-conjugacy relation between corresponding
elements of the semigroups.  More precisely, if $h=h_{j_1}\circ
\dots \circ h_{j_n} \in H$, then for $g=g_{j_1} \circ \dots \circ
g_{j_n}$ we have $k \circ h=g\circ k$ since $k\circ h_{j_{1}}\circ
\cdots \circ h_{j_{n}}
  = g_{j_{1}}\circ k\circ h_{j_{2}}\circ \cdots \circ h_{j_{n}}
  = g_{j_{1}}\circ g_{j_{2}}\circ k\circ h_{j_{3}}\circ \cdots \circ h_{j_{n}}
  = .....
  = g_{j_{1}}\circ g_{j_{2}}\circ \cdots \circ g_{j_{n}}\circ k$.

Let $z_0$ be a point in $N(H)$ and let $\Delta$ be a small open
set in $N(H)$ containing $z_0$ such that $h(\Delta)$ has spherical
diameter less than $\epsilon$ for all $h \in H$. Denoting the
Lipschitz constant (with respect to the spherical metric)
of $k$ by $C$ (see~\cite{Be}, p.~32), we see
that for any $g=g_{j_1} \circ \dots \circ g_{j_n} \in G$ the
diameter of $g(k(\Delta))=k(h_{j_1}\circ \dots \circ
h_{j_n}(\Delta))=k(h(\Delta))$ is less than $C \epsilon$.  Hence
$k(z_0)\in N(G)$ and so we conclude that $k(N(H))\subset N(G).$

Let $z_0$ be a repelling fixed point for some $h=h_{j_1}\circ
\dots \circ h_{j_n} \in H$, but which is not a critical point of
$k$. Then for $g=g_{j_1} \circ \dots \circ g_{j_n}$ we have
$g\circ k=k \circ h$ and hence $g$ has a fixed point at $k(z_0)$
with the same multiplier as that of $h$ at $z_0$ (using the chain
rule and the fact that $g=k\circ h \circ k^{-1}$ for the branch of
$k^{-1}$ which maps $k(z_0)$ to $z_0$). Hence we have shown that
the repelling fixed points of the maps in $H$, which are not any
of the finite number of critical points of $k$, map under $k$ to
repelling fixed points of maps in $G$. Since the Julia set of a
rational semigroup is a perfect set equal to the the closure of
the set of repelling fixed points of the elements of the semigroup
(see~\cite{HM1}, Theorem 3.1 and Corollary 3.1), it then follows
that $k(J(H)) \subset J(G)$.

Since $J(H)=\CC \setminus N(H)$ and $J(G)=\CC \setminus N(G)$ the
lemma now follows from the fact that $k(\CC)=\CC$.
\end{proof}

One might expect that a result similar to Lemma~\ref{semiconj}
would hold for completely invariant Julia sets, however, we
require an additional hypothesis as noted in the following lemmas.

\begin{lemma}\label{cisemiconj}
Suppose rational functions $g, h, k$ satisfy the semi-conjugacy
relation $k \circ h=g \circ k$.  If $\tilde{S}$ is completely
invariant under $g$, then $k^{-1}(\tilde{S})$ is completely
invariant under $h$. Also, if $S$ is completely invariant under
$h$ and $k^{-1}(k(S))=S$, then $k(S)$ is completely invariant
under $g$.
\end{lemma}

The proof of Lemma~\ref{cisemiconj} follows readily from the
semi-conjugacy and will therefore be omitted.


\begin{lemma}\label{cijsemiconj}
For rational semigroups $G=\langle g_j:  j \in \mathcal{I}
\rangle$ and $H=\langle h_j: j \in \mathcal{I} \rangle$ where
there exists a rational function $k$ satisfying the semi-conjugacy
relation $k \circ h_j=g_j \circ k$ for each $j \in \mathcal{I}$,
we have $k(E(H))\subset E(G)$ (and thus $W(G) \subset k(W(H))$).
If we also have that $k^{-1}(k(E(H))=E(H)$, then $k(E(H))=E(G)$
and $W(G)=k(W(H))$.
\end{lemma}
\begin{remark}
The hypothesis $k^{-1}(k(E(H)))=E(H)$ stated above would
automatically follow from the other assumptions if, in addition,
$k$ is a (branched) Galois covering. We, however, do not require
that form of the statement because one can easily check that this
hypothesis holds in the situations we consider below.
\end{remark}
\begin{proof}
Let $h=h_{j_1}\circ \dots \circ h_{j_n} \in H$ and consider the
corresponding $g=g_{j_1} \circ \dots \circ g_{j_n} \in G$.  Since
$E(G)$ is completely invariant under $g$ and  $k \circ h = g \circ
k$, Lemma~\ref{cisemiconj} shows that the closed set
$k^{-1}(E(G))$ is completely invariant under $h$. Since $h \in H$
was arbitrary, we conclude that $E(H) \subset k^{-1}(E(G))$. Thus
$k(E(H))\subset E(G)$.

Similarly one can use Lemma~\ref{cisemiconj} to show that
$k^{-1}(k(E(H))=E(H)$ implies $E(G) \subset k(E(H))$ and so $E(G)
= k(E(H))$.  When $k^{-1}(k(E(H))=E(H)$, $k$ maps $E(H)$ in a
$\deg(k)$-to-one fashion onto $k(E(H))=E(G)$.  Since $k$ is a
rational map of global degree $\deg(k)$, it must then map
$W(H)=\CC \setminus E(H)$ onto $W(G)=\CC \setminus E(G)$ (also in
a $\deg(k)$-to-one fashion).
\end{proof}

%
%
%

\begin{example}[Counterexamples to Conjectures~\ref{conj1} and~\ref{conj2}
]\label{ex1} Let $f$ be an odd rational map such that $f(U)=U$.
Then by Lemma~\ref{tildef} there exists a rational function
$\tilde{f}$ satisfying the semi-conjugacy relation $\phi \circ f =
\tilde{f} \circ \phi$. 
Similarly we let $g$ be an odd rational map with $g(U)=U$ and so
there exists a rational map $\tilde{g}$ with $\phi \circ g =
\tilde{g} \circ \phi$.  
By choosing $f$ and $g$ such that $J_f\neq \CR$ and $J_g=\CR$,
 we have
that $J_{\tilde{f}} \neq I$ and $J_{\tilde{g}}=I$ by
Lemma~\ref{semiconj}. Since $\CR$ is completely invariant under
both $f$ and $g$ we have $E(G)\subset \CR$ where $G=\langle f, g
\rangle$.  Since $E(G) \supset J_g=\CR$, we conclude that
$E(G)=\CR.$   For $\tilde{G}=\langle \tilde{f}, \tilde{g} \rangle$
we see that since
$\phi^{-1}(\phi(E(G)))=\phi^{-1}(\phi(\CR))=\CR=E(G)$, we must
have $E(\tilde{G})=\phi(\CR)=I$.  Since $J_{\tilde{f}} \neq
J_{\tilde{g}},\,\,\tilde{G}$ is a counterexample to
Conjectures~\ref{conj1} and~\ref{conj2}.

Specifically we may select $f(z)=2z-\frac{1}{z}$ and
$g(z)=\frac{z^2-1}{2z}$.  Hence $J_f$ is a Cantor subset of $I$
(see~\cite{Be}, p.21). Since $g$ is the conjugate of $z \mapsto
z^2$ under $z \mapsto i\frac{1+z}{1-z}$ we see that $J_g=\CR$.  In
this case one can calculate (via the proof of Lemma~\ref{tildef})
that $\tilde{f}(z)=\frac{3z+5z^2}{1+3z+4z^2}$ and
$\tilde{g}(z)=2z^2-1$.
\end{example}


In the next example, we construct a semigroup $G$ that provides a
counterexample to Conjectures~\ref{conj1} and~\ref{conj2} with the
additional property that $J(G) \subsetneq E(G)$.

\begin{example}
Consider $f(z)=2z - 1/z$ as in Example~\ref{ex1}. Let
$\varphi(z)=2z$, and set $g(z)=(\varphi \circ f \circ
\varphi^{-1})(z)=2z-4/z.$ Note that $J_g=\varphi(J_f)=2J_f$ and
that $\CR$ is completely invariant under $g$.  Hence for
$G=\langle f,g \rangle$, we have $E(G) \subset \CR$.

Suppose that $E(G) \neq \CR$. Since $\CR$ is completely invariant
under both $f$ and $g$, it follows from Lemma~3.2.5 in \cite{RS}
that if $E(G)$ contains a non-degenerate interval in the real
line, then $E(G)=\CR.$  Hence we may select an open interval
$L=(x,y)$ in $\R \setminus E(G)$ with both $x,y$ large. Since the
length of the intervals $f^n(L)$ tends to $+\infty$, we may assume
that $y-x$ is large.  By expanding the interval we may also assume
that $x,y \in E(G)$ (note that we used here that $\infty$ is a
non-isolated point in $E(G)$ which follows since $2 \in J_g
\subset E(G)$ and $f^n(2) \to \infty$).

Since $x$ is large, we can use the fact that $f(x)$ is slightly
greater than $g(x)$ to see that $g^{-1}(\{f(x)\})$ contains a
point slightly larger than $x$ (and hence less than $y$). But by
the complete invariance of the set $E(G)$ under $f$ and $g$, we
get $g^{-1}(\{f(x)\}) \subset E(G)$. This is a contradiction since
the interval $(x,y)$ does not meet $E(G)$.  We conclude that
$E(G)=\CR$.

Since $\infty$ is an attracting fixed point under both $f$ and
$g$, we see that small neighborhoods of $\infty$ map inside
themselves under each map in $G$.  Hence $\infty \in N(G)$ and so
$J(G) \neq \CR$.

As in Example~\ref{ex1} we may semi-conjugate the odd rational
maps $f$ and $g$ by $\phi$ to get maps
$\tilde{f}(z)=\frac{3z+5z^2}{1+3z+4z^2}$ and
$\tilde{g}(z)=\frac{5z^2+40z-29}{3z^2+40z-27}$. Hence for
$\tilde{G}=\langle \tilde{f}, \tilde{g} \rangle$ we have
$J(\tilde{G}) = \phi(J(G)) \subsetneq \phi(\CR) = I$ and
$E(\tilde{G}) = \phi(E(G)) = \phi(\CR) = I$.  Since $J_{\tilde{f}}
\neq J_{\tilde{g}}$ (otherwise one would have
$E(\tilde{G})=J(\tilde{G})=J_{\tilde{f}}=J_{\tilde{g}})$, we see
that $\tilde{G}$ is a counterexample to Conjectures~\ref{conj1}
and~\ref{conj2}.
\end{example}


\section{Counterexamples to Conjecture~\ref{Conj:Jf=Jg}}\label{counterNA}

In~\cite{HM1}, p.~366 Hinkkanen and Martin give the following
definition.
\begin{definition} \label{D:nearabel}
        A rational semigroup $G$ is {\bf nearly abelian} if there is a
        compact family of M\"obius transformations $\Phi = \{\phi\}$
        with the following properties:

        (i) $\phi(N(G)) = N(G)$ for all $\phi \in \Phi$, and

        (ii) for all $f, g \in G$ there is a $\phi \in \Phi$ such that
        $f \circ g = \phi \circ g \circ f$.
\end{definition}



\begin{theorem}[\cite{HM1}, Theorem~4.1] \label{T:nearabel1}
        Let $G$ be a nearly abelian semigroup. Then for each $g \in G$ we
        have $J_g = J(G)$.
\end{theorem}

A natural question is to what extent does the converse to
Theorem~\ref{T:nearabel1} hold.  Using a result of A.~Beardon
(see~\cite{Be1}, Theorem 1) Hinkkanen and Martin have proved the
following result for polynomial semigroups.

\begin{theorem}[\cite{HM1}, Corollary~4.1] \label{C:polysameJ}
        Let $\mathcal{F}$ be a family of polynomials of degree at
        least 2, and suppose that there is a set $J$ such that $J_g =
        J$ for all $g \in \mathcal{F}$.  Then $G = \langle \mathcal{F}
        \rangle$ is a nearly abelian semigroup.
\end{theorem}

%
%
%
%

Note that under the hypotheses of Theorem~\ref{C:polysameJ} we
have $J_h = J(G)$ for each generator $h \in \mathcal{F}$. So we
see that Conjecture~\ref{Conj:Jf=Jg} is suggesting that if
$J_h=J(G)$ for just one $h \in G$, then $G$ is still nearly
abelian.  However, this is not the case as we see by the following
counterexample.

\begin{example}[Counterexample to Conjecture~\ref{Conj:Jf=Jg}]\label{exnotNA} Let
$f(z)=z^2-2, g(z)=4z^2-2$ and $G=\langle f, g \rangle$. It is well
known that $f$ is a conjugate of $2z^2-1$ by $z \mapsto 2z$ and so
$J_f=[-2,2]$ (see~\cite{Be}, p.~9).  It can easily be seen that
$g$ maps $[-1,1]$ onto $[-2,2]$ in a two-to-one fashion.  Since
$g^{-1}([-1,1]) \subset g^{-1}([-2,2]) = [-1,1]$ it follows that
$J_g \subset [-1,1]$. In particular $J_g \subsetneq J_f$. We also
note that $\CC \setminus [-2,2]$ is forward invariant under both
$f$ and $g$ and as such must lie in $N(G)$ by Montel's Theorem. It
follows that $J(G)=[-2,2]=J_f$, yet $J_f \neq J_g$.

We remark that any map $g$ that maps a proper sub-interval of
$[-2,2]$ onto $[-2,2]$ in a $deg(g)$-to-one fashion would suffice
in the above example and such functions can easily be obtained by
constructing real polynomials with appropriate graphs.  Also, $f$
may be replaced by any Tchebycheff polynomial (see section 1.4 of
~\cite{Be}), normalized so that $J_f=[-2,2]$.
\end{example}

\textit{Acknowledgements:} The authors would like to thank Aimo
Hinkkanen for his advice given during the preparation of this
manuscript.

\bibliographystyle{plain}


\end{document}